\documentclass[12pt,twoside]{amsart}

\usepackage{verbatim}
\usepackage{amssymb,latexsym}
\usepackage{a4wide}
\usepackage[latin1]{inputenc}
\usepackage[T1]{fontenc}
\usepackage{times}

\def\M{\mathcal{M}}

\def\hpq0{h^{p,q}_{\leq 0}}
\def\Hpq0{\H_{\leq 0}^{p,q}}

\def\I{\mathcal{I}}

\def\dbar{\bar\partial}
\def\ddbar{\partial\dbar}

\def\R{{\mathbb R}}

\def\C{{\mathbb C}}

\def\F{{\mathcal F}}

\def\D{\mathcal{D}}

\def\H{{\mathcal H}}
\def\E{{\mathcal E}}

\def\Re{{\rm Re\,  }}

\def\be{\begin{equation}}
\def\ee{\end{equation}}

\newtheorem{thm}{Theorem}[section]

\newtheorem{cor}[thm]{Corollary}
\newtheorem{prop}[thm]{Proposition}

\theoremstyle{definition}

\theoremstyle{remark}

\newtheorem{preremark}{Remark}
\newtheorem{preex}{Example}

\numberwithin{equation}{section}

\title[]
{Complex Brunn-Minkowski theory and positivity of vector bundles.}

\address{Department of Mathematics\\Chalmers University
  of Technology \\
 S-412 96
  G\"OTEBORG\\SWEDEN} 

\email{ bob@chalmers.se}

\author[]{Bo Berndtsson}
\thanks{ Partially supported by grants from  Vetenskapsr\aa det}

\begin{document}

\begin{abstract}This is a survey of results on positivity of vector bundles, inspired by the Brunn-Minkowski and Pr\'ekopa theorems. Applications to complex analysis, K\"ahler geometry and algebraic geometry are also discussed. 
\end{abstract}

\maketitle

\vspace{-0.5 truecm}

\section{Introduction.}

The classical Brunn-Minkowski theorem is an inequality for the volumes of convex sets. It can be formulated as a statement about how the volumes of the vertical $n$-dimensional sections of a convex body in $\R^{n+m}$ vary with the section; more precisely it says that the $n$:th root of the volumes define a concave function on $\R^m$.  
The theorem has many applications, and it has also been generalized in many different directions ( see e. g. the survey, \cite{Gardner}). 

One important generalization is Pr\'ekopa's theorem, \cite{Prekopa}, which can be seen as a version  of the Brunn-Minkowski theorem for convex functions instead of convex sets. Let $\phi(x,t)$ be a convex function on $\R^n_x\times \R^m_t$, satisfying some mild extra conditions. We define a new function on $\R^m$ by
$$
e^{-\tilde\phi(t)} =\int_{\R^n} e^{-\phi(t,x)} dx.
  $$
 Then Pr\'ekopas theorem says that $\tilde\phi$ is also convex. Measures of the type
 $$
 e^{-\phi} dx,
 $$
where $\phi$ is convex,  are  called {\it log concave },  and in this terminology Pr\'ekopa's theorem says that the marginals, or push forwards,  of log concave measures are log concave. If we admit convex functions that attain the value $\infty$, the Brunn-Minkowski theorem is a direct consequence of Pr\'ekopa's theorem, corresponding to the case when $\phi$ is the indicator function of a convex body, i. e. the function which is zero on the convex body and infinity outside. (Properly speaking we get a version of the Brunn-Minkowski theorem with $n$:th roots replaced by logarithms. This is the 'multiplicative version' of the Brunn-Minkowski theorem, and is easily seen to be equivalent to the 'additive version'.)
 
 Among the many different proofs of Pr\'ekopa's theorem, the one that is relevant to us here is the proof of Brascamp and Lieb, \cite{Brascamp-Lieb}. It goes by first proving  the {\it Brascamp-Lieb inequality}, which is a Poincar\'e inequality for the $d$-operator (on functions). In the case when $n=1$ (which is really the main case) it says the following: Let $u$ be a smooth function on $\R$, satisfying
 $$
 \int_\R u e^{-\phi(x)} dx =0,
 $$
where $\phi$ is a smooth function, which is strictly convex in the sense that $\phi''>0$. Then
$$
\int_\R |u|^2 e^{-\phi} dx\leq \int_\R\frac{|u'|^2}{\phi''} e^{-\phi}dx.
$$
From this, Pr\'ekopa's theorem (and therefore also the Brunn-Minkowski theorem) follows just by differentiating $\tilde\phi$ twice, and applying the Brascamp-Lieb inequality to estimate the result. Actually, the two results are  `equivalent' in the sense that the Brascamp-Lieb theorem follows from Pr\'ekopa's theorem too. ( See \cite{Dario-Boaz}, which also gives a nice  complement to the oversimplified historical picture described above.)

Now we observe that the Brascamp-Lieb inequality is the real variable version of H\"ormander's $L^2$-estimate for the $\dbar$-operator,
\cite{Hormander}, \cite{Demaillynotes}. (This point of view was probably first stated clearly in \cite{Dario}.) To describe the simplest case of H\"ormander's estimate we let $\phi$ be a smooth strictly subharmonic function in $\C$. Then we let $u$ be a smooth function on $\C$, satisfying
$$
\int_\C u \bar h e^{-\phi} d\lambda=0,
$$
for all holomorphic functions $h$ satisfying the appropriate $L^2$-condition ( $d\lambda$ denotes Lebesgue measure). Note that this is a direct counterpart to the condition on $u$ in the real case.  Then $u$ was assumed to be orthogonal to all constant functions, i. e. all functions in the kernel of $d$, whereas in the complex case $u$ is assumed to be orthogonal to all functions in the kernel of $\dbar$. The conclusion of H\"ormander's estimate is now that
$$
\int_\C |u|^2 e^{-\phi} d\lambda\leq \int_\C \frac{|\dbar u|^2 }{\Delta\phi} d\lambda,
$$
which  clearly is very similar to the conclusion of the Brascamp-Lieb inequality. 
The condition that $u$ is orthogonal to all  holomorphic functions, means that $u$ is the $L^2$-minimal solution to a $\dbar$-equation, and this is how H\"ormander's theorem is mostly thought of; as an estimate for solutions to the $\dbar$-equation. (There is also an even more important existence part of H\"ormander's theorem, but that plays no role here.)
In the same way, the Brascamp-Lieb theorem is an $L^2$-estimate for the $d$-equation, and it can be obtained as a special case of H\"ormander's theorem, when the functions involved do not depend on the imaginary part of $z$.

Given the importance of Pr\'ekopa's theorem in convex geometry, it now becomes a natural question if there are any analogous consequences of H\"ormander's theorem in the complex setting, that generalize the real theory.  This is  the subject of the work that we will now describe.

The most naive generalization would be that letting $\phi$ be plurisubharmonic (i. e. subharmonic on each complex line) in $\C^n\times \C^m$, and defining $\tilde\phi$ as we did in the  real case, we get a plurisubharmonic function. This is however in general not the case as shown by a pertinent example of Kiselman, \cite{Kiselman},
$$
\phi((t,z):=|z-\bar t|^2 -|t|^2=|z|^2-2\Re zt
$$
(in $\C^2$). It turns out that instead we should think of the volume of a (convex) domain as the squared $L^2$-norm of the function 1, and the integrals of $e^{-\phi}$ as squared weigthed $L^2$-norms. Then it becomes natural to consider  $L^2$-norms of holomorphic functions in the complex case, i. e. to look at norms
$$
\|h|_t^2:= \int_{\C^n} |h(z)|^2 e^{-\phi(t,z)} d\lambda(z),
$$  
or similar expressions where we integrate over slices of pseudoconvex domains in $\C^n$ instead of the total space. Let 
$A^2_t$ be the (Bergman) space of holomorphic functions  with finite norm. Then we get a family of Hilbert spaces, indexed by $t$, i. e. a vector bundle, or at least a `vector bundle like' object (the bundles obtained are in general not locally trivial).  The theorems that we are going to discuss amount to saying that the curvature of these bundles is positive, or at least non negative, under natural assumptions. If, intuitively, we think of  the curvature as  (the negative of) the Hessian of the (logarithm of) the metric, this can be seen as a counterpart to the Brunn-Minkowski-Pr\'ekopa theorem. One can also recover Pr\'ekopa's theorem as a special case, see section \ref{main}. 

Let us add one remark on  the relation of positive curvature  to convexity in the real setting. A convex function on $\R^n$ is the same thing as a plurisubharmonic function on $\C^n$ that does not depend on the imaginary part of $z$. But, it is {\it not} the case that  
$$
-\log \|h_t\|_t^2,
$$
is plurisubharmonic in $t$ if $h_t$ is a holomorphic section of a holomorphic hermitian  vector bundle of positive curvature. This does however hold if the rank of the vector bundle is 1, so that we have a line bundle, if we also assume that $h_t\neq 0$ for all $t$.  This is precisely the situation in the real setting: The bundle of constants has rank 1 and the `section' 1 is never zero, and that is why we get simpler statements in the real setting. If we imagine a vector valued Brunn-Minkowski theory (see \cite{2Raufi}) we would get a situation similar to the complex case since 
$$
-\log \sum_1^N e^{-\phi_j}
$$
is in general not convex for convex $\phi_j$, except when $N=1$.

Looking at the complex situation, it is clear that the restriction to linear sections of a domain in Euclidean space is not as natural as in the real case. The general picture involves two complex manifolds $X$ and $Y$ of dimensions $n+m$ and $m$ respectively, and a holomorphic surjective map $p:X\to Y$ between them. This corresponds  to the previous situation when $X=\C^n\times\C^m$, $Y=U\subset\C^m$ and $p$ is the linear projection from $\C^n\times\C^m$ to $\C^m$. The role of the linear sections is  played by the fibers $X_y=p^{-1}(y)$ of the map. To get enough holomorphic objects to apply the theory to, we will also   replace holomorphic functions by  holomorphic sections of a line bundle, $L$, over $X$, and the plurisubharmonic function $\phi$ now corresponds to a hermitian metric of non negative curvature on $L$. We then have almost all the ingredients for the complex theory, but one item remains to sort out: How do we define $L^2$-norms over the fibers?

For a general complex manifold, like the fibers of the map $p$, there is of course no substitute for Lebesgue measure. The way out is to consider, instead of sections of $L$, $(n,0)$-forms on the fibers with values in $L$. Such forms with values in $L$ have natural $L^2$-norms, defined by wedge product and the metric on $L$. The bundle metric we get is
$$
\|u\|^2_y=c_n\int_{X_y} u\wedge \bar u e^{-\phi}.
$$
Again, in our model situation of Euclidean space, this corresponds  to integration with respect to Lebesgue measure. With this list of translations we obtain a  counterpart to Pr\'ekopa's theorem in the complex setting, under natural convexity assumptions on $X$. 

As it turns out, the nicest situation is when the map is proper, so that the fibers are compact, and also a submersion, so that the fibers are manifolds. In this case, and assuming also that the line bundle is trivial,  the theorem was already known. Indeed it is a part of Griffiths' monumental theory of  variations of Hodge structures, \cite{Griffiths}, \cite{Fujita}. Griffiths' point of view however was rather different. He considered the vector bundle that we are discussing, with fibers $H^{n,0}(X_y)$, as a subbundle of the Hodge bundle with fibers $H^n_{dRh}(X_y)$, with connection induced by the Gauss-Manin connection on the Hodge bundle. It seems difficult to generalize this approach completely to the twisted case,  when $L$ is nontrivial, since e.g.  there is no twisted counterpart to the Hodge bundle (see however \cite{Kawamata}). It is also interesting to notice the different roles played by holomorphic forms: Griffiths was interested in holomorphic forms per se because of their relation to the period map. Here the forms are forced upon us in order to define $L^2$-norms.

In the next section we will give more precise formulations of the basic results. After that we will turn to applications; to $L^2$-extension problems for holomorphic functions, the (Mabuchi) space of K\"ahler metrics in a fixed cohomology class, variations of complex strucures, and finally positivity of direct image bundles in algebraic geometry.

\bigskip

{\bf Acknowledgement:} The work described here has been influenced by discussions with many people. In particular I owe a great debt to Robert Berman, Sebastien Boucksom, Dario Cordero-Erausquin, Bo'az Klartag, Mihai Paun, Yanir Rubinstein and Xu Wang. 

\section{ The basic results}\label{main}

We begin with non proper maps, and start by looking at the simplest case. We let $\Omega$ be a  domain in $\C^n$ and $U$ a domain in $\C^m$. We then let $X=\Omega\times U$ and let $p:X\to U$ be the linear projection on the second factor. We will assume that $\Omega$ is pseudoconvex (meaning that it has a smooth strictly plurisubharmonic exhaustion function). Let $\phi=\phi(z,t)$ be  plurisubharmonic and smooth up to the boundary. Let
$$
A^2(\Omega, \phi)=\{u\in H(\Omega); \int_\Omega |u|^2 e^{-\phi(\cdot, t)}  < \infty \},
$$
be the corresponding Bergman space of holomorphic functions, equipped with the norm
$$
\|u\|^2_t=\int_\Omega |u|^2 e^{-\phi(\cdot, t)}.
$$
Since  $\phi$ is bounded, all the Bergman spaces are the same as linear spaces, but their norms vary with $t$. We therefore get a trivial vector bundle $E=A^2\times U$ over $U$ with a non trivial Hermitian metric. We define its complex structure by saying that a section $t\to h(t, \cdot)$ is holomorphic if $h$ is holomorphic on $\Omega\times U$. 
\begin{thm}\label{Ann}(\cite{1Berndtsson})
  The bundle $(E, \|\cdot\|_t)$ has non negative (Chern)  curvature in the sense of Nakano.
\end{thm}
There are two main notions of  positivity for Hermitian vector bundles; positivity in the sense of Griffiths and in the sense of Nakano. Positivity in the sense of Griffiths is defined in terms of the curvature tensor in the following way. Recall that the curvature is a $(1,1)$-form, $\Theta$, with values in the space of endomorphisms of $E$. Then, if $u_t\in E_t$,
$$
\langle \Theta u_t, u_t\rangle_t
$$ 
is a scalar valued $(1,1)$-form. Positivity  in the sense of Griffiths means that it is a positive form for any $u_t$. For the definition of positivity in the sense of Nakano, which is a stronger notion, we refer to \cite{Demaillynotes} or\cite{1Berndtsson}. We should also point out that we are here in the somewhat non standard situation of infinite rank bundles; this is also discussed in \cite{1Berndtsson}.

We first discuss the proof very briefly. Let $F$ be the bundle $(L^2(\Omega,\phi), \|\cdot\|_t)$, whose fibers are general $L^2$-functions, not necessarily holomorphic. It is also a trivial bundle, and we define a section $t\to u_t$ to be  holomorphic if the dependence on $t$ is holomorphic. Then $E$ is a subbundle of $F$. It is easy to see that the curvature of $F$ is the $(1,1)$-form $\sum_1^m \phi_{t_j,\bar t_k} dt_j\wedge d\bar t_k$, where the coefficients should be interpreted as the endomorphisms
$$
u_t\to \phi_{t_j, \bar t_k} u_t.
$$
Then it is immediately clear that $F$ has non negative curvature as soon as $\phi$ is plurisubharmonic in $t$ for $z$ fixed. By general principles (\cite{Demaillynotes}), the curvature of the subbundle $E$ is given by
$$
\langle \Theta^E u_t, u_t\rangle_t =\langle \Theta^F u_t, u_t\rangle_t -\|\pi_{E^\perp}\theta^F u_t\|_t^2.
$$
Here $\theta^F$ is the connection form for the Chern connection of $F$, and $\pi_{E^\perp}$ is the orthogonal projection on the orthogonal complement of $E$ in $F$. The important thing to notice is now that, since $\pi_{E^\perp}\theta^F u_t$  lies in the orthogonal complement of the space of holomorphic functions, it is the $L^2$-minimal solution of some $\dbar$-equation, and we can apply H\"ormander's estimate (for pseudoconvex domains in $\C^n$). This allows us to control the second, negative, term by the first, positive,  one, and that gives the theorem.

It is well known that positivity in the sense of Griffiths, is equivalent to negativity of the dual bundle, $E^*$. On the other hand, negativity in the sense of Griffiths, is equivalent to saying that
$$
\log \|v_t\|_t
$$
is plurisubharmonic for any holomorphic section of $E^*$.

This leads to  a more concrete reformulation of the first theorem. For $t$ in $U$, let $t\to \mu_t$ be a family of complex measures on $\Omega$. Assume there is a compact subset, $K$,  of $\Omega$, such that all $\mu_t$ are supported 
on $K$. Then
$$
u_t\to \int_\Omega u_td\mu_t=:\mu_t(u_t)
$$
defines a secion of the dual bundle $E^*$ and it gets a norm inherited from the norm on $E_t$.
\begin{cor}\label{Baby} Assume that the section $\mu_t$ is holomorphic in the sense that
  $$
  t\to \mu_t(h(t,\cdot)
  $$
  is holomorphic for any $h$ holomorphic on $X$. Then
  $$
  \log\|\mu_t\|_t
  $$
  is plurisubharmonic.
\end{cor}

But this statement makes sense in a much more  general situation. 
\begin{thm}\label{Gent}  Let $D$ be a pseudoconvex domain in $\C^n\times\C^m$, and denote by $D_t=\{z\in\C^n; (z,t)\in D\}$ its vertical slices. Let $\mu_t$ be a family of measures on $D_t$ that are all locally supported in a compact subset of $D$, and depend on $t$ in a holomorphic way  so that
$$
t\to \mu_t(h(t,\cdot))
$$
is holomorphic in $t$ if $h$ is holomorphic  in $D$. Then
$$
\log\|\mu_t\|_t
$$
is plurisubharmonic in $t$.
\end{thm}

This way we have implicitly  defined positivity of curvature for the 'bundle of Bergman spaces' $A^2(D_t,\phi)$, even though this is not properly speaking a vector bundle, since it is not locally trivial.
Here we also note that Theorems \ref{Gent} and \ref{Baby}  imply Pr\'ekopa's theorem: Take $D=(\C^*)^n$ (where $\C^*=\C\setminus\{0\}$), and define $\mu$ by taking averages over the $n$-dimensional real torus
$$
\mu(h)=\int_{T^n} h(e^{i\theta_1}z_1,..e^{i\theta_n}z_n)d\theta
$$
(it does not depend on $z$). 
Computing the norm of $\mu$  as a functional on
$A^2(D,e^{-\phi(t,\cdot)})$, where $\phi$ only depends on $|z_j|$,   we recover Pr\'ekopa's theorem.

One main case of Theorem \ref{Gent} is when the measures $\mu_t$ are all Dirac delta functions.
\begin{thm}\label{AIF} (\cite{2Berndtsson}) Let $D$ be a pseudoconvex domain in $\C^n\times\C^m$ and let $\phi$ be plurisubharmonic in $D$. For each $t$ in the projection of $D$ to $\C^m$ let $B_t(z)=B_t(z,z)$ be the diagonal Bergman kernel for $A^2(D_t, \phi)$. Then $\log B_t(z)$ is plurisubharmonic in $D$.
\end{thm}
This theorem was obtained earlier when $n=1$ and $\phi=0$ by Maitani and Yamaguchi in \cite{Maitani-Yamaguchi}. It is a  consequence of the previous result. Let $t\to f(t)$ be a holomorphic map such that $f(t)\in D_t$, and let
$\mu_t$ be a Dirac mass at $f(t)$. It is immediate that $t\to \mu_t(h)=h(t,f(t))$ is holomorphic if $h$ is holomorphic on the product space.  Moreover, 
$$
\|\mu_t\|^2_t= B_t(f(t)),
$$
so by the theorem,  $\log B_t(f(t))$ is plurisubharmonic for any such map. This  means  that $B_t(z)$ is plurisubharmonic in $D$.

The proof of Theorem \ref{Gent} uses the corollary, applied to a sequence of weights that tend to infinity outside of $D$. It is given in detail in \cite{2Berndtsson} for the situation of Bergman kernels, but the same proof gives the more general result too. See also \cite{1Wang} for a detailed analysis of a curvature operator defined in the setting of Theorem \ref{Gent} (under one extra assumption) and conditions for when the curvature vanishes.

We next turn to general surjective holomorphic maps between complex manifolds, and there we will restrict to the case of proper maps, to fix ideas. Let $X$ be a complex manifold of dimension $n+m$ and $U$ a complex manifold of dimension $m$. Since our results are mostly local we can here think of $U$ as a domain in $\C^m$. Let $p:X\to U$ be a holomorphic and surjective map. We say that $p$ is smooth if its differential is surjective at every point. Then the fibers $X_t=p^{-1}(t)$ are smooth manifolds, and when the map is proper they are also compact. Let $L\to X$ be a holomorphic line bundle over $X$, equipped with a an Hermitian metric $e^{-\phi}$.

For $t$ in the base, $U$, we let
$$
E_t=H^{n,0}(X_t, L),
$$
be the space of holomorphic $L$-valued $(n,0)$-forms on $X_t$. 
We can also think of this space as
$$
E_t=H^0(X_t, K_{X_t}+L),
$$
the space of holomorphic sections over $X_t$ of the canonical bundle of $X_t$ twisted with the bundle $L$.  Equivalently
$$
E_t=H^0(X_t, K_{X/U}+L),
$$
is the space of  holomorphic sections over $X_t$ of the relative canonical bundle $K_{X/U}$ twisted with $L$. (The relative canonical bundle, defined as $K_{X/U}= K_X-K_U$, is a line bundle over the total space $X$ that restricts to $K_{X_t}$ on every fiber $X_t$.) Let
$$
E=\cup_{t\in U} \{t\}\times E_t.
$$
It turns out that, when the metric on $L$ has semipositive curvature, so that $i\ddbar\phi\geq 0$, then $E$ is a holomorphic vector bundle over $U$. A section $u_t$  of $E$ is defined to be holomorphic if
$$
u_t\wedge dt
$$
is a holomorphic $(n+m,0)$ -form on $X$, where $(t_1, ...t_m)$ are holomorphic coordinates on $U$ and $dt=dt_1\wedge ... dt_m$. Finally, the $L^2$-norm 
$$
\|u_t\|^2_t:= c_n\int_{X_t} u_t\wedge \bar u_t e^{-\phi},
$$
defines an Hermitian structure   on $E$. We are then in a situation analogous to the one in Theorem \ref{Ann}, with the advantage that we now get a bona fide vector bundle of finite rank.
\begin{thm}\label{Ann2}(\cite{1Berndtsson}) Assume  that the curvature  of $L$ is nonnegative ( $i\ddbar\phi\geq 0$) and that the total space is K\"ahler. Then  $E$ is a holomorphic vector bundle with nonnegative curvature in the sense of Nakano (and hence also in the sense of Griffiths).
\end{thm}
Notice that what corresponds to pseudoconvexity here is the assumption that $X$ is K\"ahler, so somewhat curiously the K\"ahler assumption appears as a convexity condition. If we look at the global picture and assume that $X$ is even projective this is rather natural, since we may then remove a divisor from $X$ and get a total space that is Stein, hence pseudoconvex. Since a divisor is removable for $L^2$-holomorphic forms, it appears that projectivity is related to (pseudo)convexity for these problems, but as it turns out, K\"ahler is enough.

When the fibration is trivial, so that $X=Z\times U$, the theorem can be proved in much the same way as Theorem \ref{Ann}, at least when the curvature of $e^{-\phi}$ is strictly positive on fibers. The proof in the general case is based on the formalism of pushforwards of currents. We write the norm squared of a holomorphic section as the pushforward of a form on the total space under $p$, compute $i\ddbar\|u_t\|^2_t$ using this formalism, and can then identify the curvature.

Theorem \ref{Ann2} can be developed in two different ways, by  either adding more, or assuming less,  assumptions. Let us first assume that the curvature form of the metric on $L$, $i\ddbar\phi$ restricts to a strictly positive form on each fiber. Then we have a K\"ahler metric $\omega_t:=i\ddbar\phi|_{X_t}$ on each fiber, and thus an associated Laplace operator, $\Box_t=\dbar\dbar^*+\dbar^*\dbar$  on $L$-valued forms
on each fiber. We can then give an 'explicit' formula for the curvature. To formulate the result we also assume that the base dimension $m=1$; the general case follows by restricting to lines or curves in the base. We then define a $(1,1)$-form
$$
C(\phi)= c(\phi) idt\wedge d\bar t,
$$
by
$$
(i\ddbar\phi)^{n+1}/(n+1)!=C(\phi)\wedge (i\ddbar\phi)^n/n!
$$
(we will come back to this form later in sections 4 and 5).
\begin{thm}\label{5pos}(\cite{3Berndtsson})
  Assume that $\omega_t=i\ddbar\phi|{X_t} >0$ on each fiber $X_t$. Then the curvature, $\Theta$ of $(E, \|\cdot\|_t)$ is given by
  $$
  \langle \Theta_{\partial/\partial t,\partial/\partial \bar t} \, u_t, u_t\rangle_t= c_n\int_{X_t} c(\phi) u_t\wedge\bar u_t e^{-\phi} +\langle (1+\Box_t)^{-1} \kappa_\phi\cup  u_t,\kappa_\phi\cup u_t\rangle_t
  $$
  (where we define $\kappa_\phi\cup$ below).
\end{thm}
To define $\kappa_\phi$, we first have to define the notion of 'horizontal lift' of a vector field, $\partial/\partial t$ on the base, as introduced by Schumacher, \cite{Schumacher}, generalizing the earlier notion of 'harmonic lift' of Siu, \cite{Siu}. A vector field on the total space $X$ is said to be vertical, if it maps to zero under $p$. It is horizontal if it is orthogonal to all vertical fields under  the (possibly degenerate, possible indefinite) metric $i\ddbar\phi$. As shown in \cite{Schumacher}, any vector field on the base $U$ has a unique horizontal lift. In our case, the vector field on the base is $\partial/\partial t$, and there is a unique horisontal field $V_\phi$ (depending on $i\ddbar\phi$) which maps to $\partial/\partial t$ under $p$. Taking $\dbar V_\phi$ and restricting to fibers we get $\kappa_\phi$; it is a $\dbar$-closed $(0,1)$-form on each fiber, with values in $T^{1,0}(X_t)$. The cohomology class of $\kappa_{\phi}$ in $H^{0,1}(X_t, T^{1,0}(X_t))$ is the Kodaira-Spencer class, but here it is important to look at this particular representative of the class. The cup product $\kappa_\phi\cup u_t$ is obtained by contracting with the vector part of $\kappa_\phi$ and wedging with the form part. 

From the formula we see that if $i\ddbar\phi\geq 0$ and the curvature is zero, then $c(\phi)=0$ and $\kappa_\phi$ is zero on each fiber. These two conditions imply that the horizontal lift of $\partial/\partial t$ is a holomorphic field, whose flow maps fibers to fibers. In fact, we have
\begin{cor}\label{5pos2} Assume that $i\ddbar\phi\geq 0$, $i\ddbar\phi >0$ on fibers and that the curvature $\Theta$ vanishes. Then there is a holomorphic vector field $V$ on $X$ whose flow lifts to $L$, such that its flow maps fibers of $p$ to fibers and is an isometry on $L$. In particular, the flow maps $\omega_t$ to $\omega_{t'}$.
  \end{cor}  
In general terms, this means that the curvature can only vanish if the fibration $p:X\to U$ is trivial, holomorphically and metrically. This is proved under the assumption that the curvature of $L$ is strictly positive on each fiber, and it does not hold in general without this assumtion (see \cite{3Berndtsson}).

It is also interesting to compare again to Pr\'ekopa's theorem. There, if the function $\tilde\phi$ is not strictly convex but linear (and if we assume that $m=1$), then $\phi$ must have the form
$$
\phi(t,x)=\psi(x+tv)+ct,
$$
where $\psi$ is a convex function on $\R^n$,  $v$ is a  vector in $\R^n$ and $c$ is a constant, \cite{Dubuc}. Thus, the variation of $\phi$ with respect to $t$ comes from the flow of a constant vector field applied to a fixed function $\psi$, 'lifted' to the line bundle $\R^n\times \R$ by adding $ct$. This is similar to what happens in the complex setting, except that there we get a holomorphic vector field instead.

We next discuss versions of Theorem \ref{Ann2} in the more general setting when the metric on the line bundle $L$ is allowed to be singular, and the map $p$ is no longer assumed to be smooth (in the sense described above) but only surjective. On the other hand, we now assume that $X$ is projective and that $p:X \to Y$, where $Y$ is a compact manifold.

By Sard's theorem, $p$ is smooth outside of $p^{-1}(Y_1)$, where $Y_1$ is a proper analytic subvariety  of $Y$, and there we can define the $L^2$-metric as before. 
There is a counterpart of the Bergman kernel here. It is first defined for $a\in p^{-1}(Y_1)$ as the norm of the  evaluation functional at $a$ on $E_{p(a)}$,
$$
B(a)=\sup |u(a)|^2/\|u\|^2_{p(a)},
$$
where the supremum is taken over all sections $u$ in $H^0(X_{p(a)}, K_{X/Y}+L)$ that are square integrable with respect to the metric $e^{-\phi}$ (if there are no nontrivial such section, we let $B(a)=0$). 
This definition however depends on the trivialization of $K_{X/Y}+L$ chosen near $a$, so the Bergman kernel is not a function but defines a metric on $K_{X/Y}+L$ -- the Bergman kernel metric $B^{-1}=e^{-\log B}$. 
\begin{thm}\label{Duke} (\cite{Berndtsson-Paun})If the singular metric $e^{-\phi}$ has semipositive curvature (i. e. $i\ddbar\phi\geq 0$ in the sense of currents), and the Bergman kernel metric is not identically equal to 0, it defines a singular metric of semipositive curvature on $p^{-1}(Y_1)$. Moreover, this metric extends to a singular metric on  $K_{X/Y}+L$ over all of $X$.
  \end{thm}
Briefly, if $L$ is pseudoeffective, i. e. has a singular metric of nonnegative curvature, $K_{X/Y}+L$ is also pseudoeffective, provided that it has a non trivial $L^2$-section over at least one fiber.

This result can be seen as a counterpart to Theorem \ref{Gent} in this setting. The more difficult problem of counterparts of Theorem \ref{Ann2} is discussed in the last section. 
\section{ The Suita conjecture and $L^2$-extension.}

Let $D$ be a (say smoothly) bounded domain in $\C$ and suppose $0\in D$. The Bergman space of $D$ is 
$$
A^2(D)=\{h\in H(D); \|h\|^2:=\int_D |h|^2 d\lambda<\infty\},
$$
and the Bergman kernel at 0 is
$$
B(0)=\sup_{\|h\|\leq 1} |h(0)|^2.
$$
To state Suita's conjecture we also need the Green's function $G(z)$ which is a subharmonic function in $D$, vanishing at the boundary, with a logarithmic pole at 0;
$$
G(z)=\log|z|^2 -v(z),
$$
where $v$ is harmonic and chosen so that $G$ vanishes at the boundary. By definition, $v(0):=c_D(0)=c_D$ is the Robin constant of $D$ at 0. Suita's conjecture, which was proved in \cite{Blocki}, \cite{Guan-Zhou}, says that
$$
B(0)\geq \frac{e^{-c_D}}{\pi}.
$$
When $D$ is a disk with center 0 it is clear that equality holds, and it was proved in \cite{Guan-Zhou} that equality holds only then. In \cite{2Blocki}, B{\l}ocki gave a much simpler proof of the conjecture, based on variations of domains and the 'tensor power trick'. In connection with this, L\'aszl\'o Lempert proposed an even simpler proof, using (pluri)subharmonic variation of Bergman kernels from \cite{Maitani-Yamaguchi}, \cite{2Berndtsson}, that we shall now describe (see \cite{2Blocki}, \cite{Berndtsson-Lempert}).

Let for $t\leq 0$,
$$
D_t=\{z\in D; G(z)<t\},
$$
and $B_t(0)$ be the Bergman kernel for $D_t$. Since
$$
\D:=\{(\tau,z)\in \C\times D; G(z) -\Re \tau<0\}
$$
is a pseudoconvex domain in $\C^2$, and $D_t=\D_\tau$, the vertical slice of $\D$,  if $t=\Re\tau$, it follows from \cite{Maitani-Yamaguchi}, \cite{2Berndtsson} (cf Theorem \ref{AIF}) that $\log B_t(0)$ is convex. When $t$ approaches $-\infty$, $D_t$ is very close to a disk centered at the origin with radius $e^{(t+c_D)/2}$, so
$$
B_t\sim \frac{e^{-t-c_D}}{\pi}.
$$
In particular, the function
$$
k(t):= \log B_t(0)+t
$$
is convex and bounded on the negative half axis.  Therefore, it must be increasing. Hence
$$
B(0)=e^{k(0)}\geq \lim_{t\to -\infty} e^{k(t)}=\frac{e^{-c_D}}{\pi},
  $$
which proves Suita's conjecture.

It is clear that similar proofs work in higher dimensions and also for weighted Bergman spaces with a plurisubharmonic weight function. The main new observation in \cite{Berndtsson-Lempert} is that the same technique can be used to give a proof of the Ohsawa-Takegoshi extension theorem, \cite{Ohsawa-Takegoshi}.

We first recall a simple version of this important result. We  let $D$ be a bounded pseudoconvex domain in $\C^n$ and $\phi$ a plurisubharmonic function in $D$. Let $V$ be a linear subspace of $\C^n$ of codimension $m$ which intersects $D$. The Ohsawa-Takegoshi theorem says that in this situation (and actually under much more general conditions), for any holomorphic function $f$ on $V\cap D$, there is a holomorphic function $F$ in $D$, such that $f=F$ on $V$ and
$$
\int_D |F|^2 e^{-\phi} d\lambda_n \leq C\int_{V\cap D} |f|^2 e^{-\phi} d\lambda_{n-m},
$$
where $C$ is a constant depending only on the diameter of $D$. What makes this theorem so powerful is that even though there is no assumption of strict plurisubharmonicity or smoothness, we get an estimate with an absolute constant. Note that, when $V$ is just a point, we get the existence of a function $F$ in $D$ with good $L^2$-estimates and prescribed value at the origin; this is equivalent to an estimate for the Bergman kernel.

Just as in Suita's conjecture, one can now ask for optimal estimates. Such estimates were given in \cite{Blocki} and more generally  by Guan and Zhou, \cite{Guan-Zhou}. We shall now see that results of this type also follow from Theorem \ref{Gent}. We also point out that, conversely, Guan and Zhou show that results along the line of Theorems \ref{Gent} and \ref{Ann2} follow from sharp versions of the Ohsawa-Takegoshi extension theorem, so in very general terms the two types of results are perhaps 'equivalent'.

Write for $z$ in $\C^n$, $z=(z_1, ...z_m, z_{m+1}, ...z_n)=(z',z'')$, and say that $V$ is defined by the equation $z'=0$. Let
$$
G(z)=\log|z'|^2,
$$
and let
$$
D_t=\{z\in D; G(z)<t\}.
$$
We may assume that $|z'|\leq 1$ in $D$. For each $t\in (-\infty, 0)$ we let $F_t$ be the extension of $f$ to $D_t$ of minimal norm. Such a minimal element exists by general Hilberts space theory. The following result is not explicitly stated in \cite{Berndtsson-Lempert}, but implicitly contained there.
\begin{prop}\label{BL}
  $$
  e^{-mt}\|F_t\|^2_t=e^{-mt}\int_{D_t} |F_t|^2 e^{-\phi} d\lambda
  $$
  is a decreasing function of $t$.
\end{prop}
It follows that
$$
\|F_0\|^2\leq \lim_{t\to-\infty} e^{-mt}\|F_t\|^2_t.
$$
It is not hard to see that the last limit equals (with $\sigma_m$ the volume of the $m$-dimensional unit ball)
$$
\sigma_m \int_{V\cap D} |f|^2 e^{\phi}d\lambda_{n-m},
$$
at least  when $\phi$ is smooth in a neighbourhood of the closure of $D$. This gives the sharp extension theorem in this setting, and the general case is  obtained by the usual approximation procedures.

Proposition \ref{BL} can be proved much like the Suita conjecture, but using the  general Theorem \ref{Gent} instead of plurisubharmonic variation of Bergman kernels.
First note that since any function $f$ on $V\cap D$ extends to $D$, the space of holomorphic functions on $V\cap D$ can be identified with the quotient space $H(D)/J(V)$, or $H(D_t)/J(V)$, where $J(V)$ is the subspace of functions that vanish on $V$. The quantity that we want to estimate, $\|F_t\|_t$, is the norm of $f$ in $A^2(D_t)/J(V)$. The space of measures with compact support in $V\cap D$ is dense in the dual of $H(D_t)/J(V)$.
We first prove a dual estimate on such measures:
$$
\|\mu\|^2_t e^{mt}
$$
is increasing, with
$$
\|\mu\|^2_t=\sup_F \frac{|\int F d\mu|^2}{\int_{D_t} |F|^2 e^{-\phi} d\lambda}.
$$
This is proved in almost exactly the same way as in the Suita case, when $V$ is a point, and Proposition \ref{BL} follows.

\section{ The space of K\"ahler metrics.}\label{Kahler}

In this section we will look at fibrations with $X=Z\times U$ where $Z$ is a compact complex manifold, $U$ is a domain in $\C$ and $p:X\to U$ is the projecion on the second factor. We also assume given a complex line bundle, $L$ over $Z$. Pulling it back to $X$ by the projection on the first factor we get a line bundle over $X$, that we denote by the same letter, sometimes writing $L\to X$ or $L\to Z$, to make the meaning clear. Mostly we assume that $L\to Z$ is positive, in the sense that it has some smooth metric of positive curvature, but we don't fix any metric -- instead we will use this set up to study the space $M_L$ of positively curved metrics on $L$.

Here all the fibers $X_t=Z$ are the same and the line bundles $L|_{X_t}$ are also the same. A metric on $L\to X$ can be interpreted as a complex curve of metrics on $L$, parametrized by $\tau$ in $U$. Of particular interest for us is the case when $U=\{\tau\in \C; 0<\Re \tau <1\}$ is a strip and $\phi$ only depends on the real part of $\tau$; then we can interpret $\phi$ as a real curve in $M_L$.

The space $M_L$ was introduced by Mabuchi, \cite{Mabuchi}. It is an open subset of an affine space modeled on $C^{\infty}(Z)$, therefore the tangent space of $M_L$ at any point $\phi$ is naturally identified with $C^{\infty}(Z)$. Mabuchi defined a Riemannian structure on $M_L$ by
$$
\|\chi\|^2_\phi:=\int_Z |\chi|^2 \omega_\phi^n,
$$
for $\chi\in C^{\infty}(Z)$, with $\omega_\phi:=i\ddbar\phi$, see also Semmes, \cite{Semmes} and Donaldson, \cite{Donaldson}. As proved in these papers, a  curve in $M_L$, i. e. a metric on $L$ depending only on $t=\Re \tau$ is a geodesic for the Riemannian metric if and only if $i\ddbar\phi\geq 0$ on $X$, and  it satisfies the homogenous complex Monge-Amp\`ere equation
$$
(i\ddbar\phi)^{n+1}=0.
$$
In our terminology in section \ref{main}, this means that $c(\phi)=0$, and in general, $c(\phi)$ can be interpreted as the geodesic curvature of the curve. Since honest geodesics are scarce, see \cite{Lempert-Vivas}, \cite{Darvas-Lempert}, \cite{David-Ross}, \cite{2David-Ross} for problems encountered in solving the homogenous complex Monge-Amp\`ere equation involved, we will also use 'generalized geodesics'. These are metrics on $L$ as above that are only assumed to be locally bounded, such that $i\ddbar\phi\geq 0$ and
$(i\ddbar\phi)^{n+1}=0$ in the sense of pluipotential theory (\cite{Bedford-Taylor}). It is easy to see (\cite{3Berndtsson})  that any two metrics in $M_L$ can be connected with a generalized geodesic in this sense, and by a famous result  of Chen, \cite{Chen}, \cite{3Blocki}, the geodesic has in fact higher regularity, so that $\ddbar\phi$ is a current with bounded coefficients. We will say that such geodesics are of class $C^{(1,1)}$. There is by now an extensive theory on these matters, for which we refer to  \cite{Mabuchi}, \cite{Semmes},  \cite{Donaldson}, \cite{Phong-Sturm} and \cite{2Donaldson}.

We first consider the case when $L=-K_Z$ is the anticanonical bundle of $Z$. The positivity of $L$ then means that $Z$ is a Fano manifold. The fibers of the vector bundle $E$, defined in section \ref{main}, 
$$
E_\tau= H^0(X_\tau, K_{X_\tau}+L)= H^0(Z,\C)
$$
are then just the space of constant functions on $Z$ (so $E$ is a line bundle). A metric $e^{-\phi}$ on $L\to Z$  can be identified with a volume form on $Z$ that we also denote $e^{-\phi}$, and with this somewhat abusive notation we have for the element '1' in $H^0(Z,\C)$, 
$$
\|1\|^2_\phi=\int_Z e^{-\phi}.
$$
We are therefore in the situation described in the introduction when $E$ is a line bundle with a nonvanishing section. Therefore we get a Pr\'ekopa type theorem for Fano manifolds.  
\begin{thm}\label{Inv}(\cite{4Berndtsson}) Let $e^{-\phi}$ be a locally bounded metric on $L\to Z\times U$, where $L\to Z$ is the anticanonical bundle of $Z$. Assume that $U$ is a strip and that $\phi$ depends only on the real part of $\tau\in U$. Suppose also that $i\ddbar\phi\geq 0$ in the sense of currents.  Let
  $$
  \tilde\phi(t)=-\log\int_Z e^{-\phi(t, \cdot)}.
  $$
  Then $\tilde\phi$ is a convex function of $t$. If $\tilde\phi$ is linear, then there is a holomorphic vector field $V$ on $U$, which is a lift of $\partial/\partial \tau$ on $U$, whose flow maps fibers $X_\tau$ to fibers and $i\ddbar\phi|_{X_\tau}$ to $i\ddbar\phi|_{X_{\tau }}$. 
 \end{thm}
When $\phi$ is smooth the first part of the theorem follows immediately from Theorem \ref{Ann2}, and the general case is easily obtained by approximation. The second part is clearly of the same vein as Corollary \ref{5pos2}, but it does not follow directly from there because of the lack of smoothness and strict positivity.

The main interest of this result lies in  its connection with K\"ahler-Einstein metrics. To explain this we also need to introduce the Monge-Amp\`ere energy, $\E(\phi)$, of a metric. This is a real valued function on $M_L$, which can be defined (up to a constant) by the property that for any curve $t\to\phi_t\in M_L$,
$$
\frac{d}{dt}\E(\phi_t)=\frac{1}{Vol(L)}\int_Z \dot\phi_t \omega_{\phi_t}^n/n!,
$$
where $\dot\phi=d\phi/dt$ and the volume of the line bundle is
$$
Vol(L):=\int_Z\omega_\phi^n/n!
$$
(it does not depend on the metric). One then defines the Ding functional, \cite{Ding},  as
$$
D(\phi):=\log\int_Z e^{-\phi} +\E(\phi).
$$
The critical points of the Ding functional are the metrics that satisfy
$$
e^{-\phi} = C\omega^n_\phi.
$$
This means precisely that the Ricci curvature of the K\"ahler metric $\omega_\phi$ equals $\omega_\phi$, i. e. that $\omega_\phi$ is a K\"ahler-Einstein metric (of positive curvature). It is well known that the Monge-Amp\`ere energy is linear along (even generalized) geodesics, so the Ding functional is linear along a geodesic precisely when its first term
$$
t\to \log\int_Z e^{-\phi_t}
$$
is linear. As noted by Berman, \cite{Berman}, this gives a proof of the Bando-Mabuchi uniqueness theorem (\cite{Bando-Mabuchi}), for K\"ahler-Einstein metrics on Fano manifolds. Indeed, suppose $e^{-\phi_0}$ and $e^{-\phi_1}$ are metrics on the anticanonical bundle of $Z$, and that $i\ddbar\phi_0$ and $i\ddbar\phi_1$ are both K\"ahler-Einstein.  Connect them by a generalized geodesic. Since both endpoints are critical points for the Ding functional, and the Ding functional is concave along the geodesic, it must in fact be linear. Then Theorem \ref{Inv} implies that $\omega_{\phi_0}$ and $\omega_{\phi_1}$ are connected via the flow of a holomorphic vector field, i. e. an element in the identity component of the automorphism group. This is the Bando-Mabuchi theorem.

As also noted by Berman, in \cite{Berman}, similar arguments prove a variant of the  Moser-Trudinger-Onofri inequality (\cite{Moser}, \cite{Onofri}) for K\"ahler-Einstein  Fano manifolds (the original case of the theorem was for $Z$ equal to the Riemann sphere). The (variant of) the theorem says, in our terminology,  that K\"ahler-Einstein metrics are global maxima for the Ding functional on the space of all positively curved metrics, which is clear from the concavity, since K\"ahler-Einstein metrics are critical points. The original version of the theorem, for the Riemann sphere, deals with metrics that are not necessarily positively curved, but as shown by Berman, this follows in one dimension from the positively curved case by an elegant trick.

The formalism described here can easily be generalized. We fix a metric $e^{-\psi}$ on a pseudoeffective $\R$-line bundle $L'$ and study metrics $e^{-\phi}$ on an $\R$-line bundle $L''$, such that $L'+L''=L=-K_Z$. With this we define the twisted Ding functional
$$
D_\psi(\phi):= \log\int_Z e^{-\phi-\psi} + \E(\phi),
$$
where Vol$(L)$ is replaced by Vol$(L'')$ in the definition of $\E$. The critical points of $D_\psi$ satisfy
$$
e^{-\phi-\psi}= C\omega_\phi^n.
$$
Hence the Ricci curvature of $\omega_\phi$ satisfies
$$
Ric(\omega_\phi)=\omega_\phi +i\ddbar\psi,
$$
so $\omega_\phi$ is a 'twisted K\"ahler-Einstein metric'. It is clear that if $\psi$ is assumed to be bounded, the same argument as before gives uniqueness modulo automorphisms, but actually the argument can be generalized to when $i\ddbar\psi=\beta [\Delta]$ is a multiple of a current defined by a divisor in $Z$. This leads to a uniqueness theorem for K\"ahler-Einstein metrics with conical singularities, introduced by Donaldson in \cite{3Donaldson}, as a tool for   the solution of the Yau-Tian-Donaldson conjecture, see \cite{CDS} and subsequent papers. In these papers it was also shown that this generalized version of Theorem \ref{Inv} can be used to prove reductivity of the group of automorphisms of the pair $(Z,\Delta)$. Furthermore, the theorem was generalized to 'log-Fano' manifolds, that arise naturally in this context,  in \cite{Berman et al}. In this context, see also \cite{2Berman} for applications to the converse direction of the Yau-Tian-Donaldson conjecture.

So far we have discussed only the case when $L$ is the anticanonical bundle of $Z$, and the resulting convexity preperties of the Ding functional, but it turns out that the formalism is also useful in connection with other functionals in K\"ahler geometry. A case in point is the Mabuchi K-energy, \cite{2Mabuchi}, $\M$. The K-energy of a metric $e^{-\phi}$ on a positive line bundle $L$ can be defined (again up to a constant) by
$$
\frac{d}{dt}\M(\phi_t)=\int_Z\dot\phi(S_\phi-\hat S_\phi) \omega_\phi^n/n!,
$$
where $t\to\phi_t$ is any smooth curve in $M_L$, $S_\phi$ is the scalar curvature of the K\"ahler metric $\omega_\phi$, and $\hat S_\phi$ is the average of $S_\phi$ over $Z$. The raison d'\^etre of $\M$ is that its critical points are precisely (the potentials of)  the metrics of constant scalar curvature.
It was proved by Mabuchi \cite{2Mabuchi}, that $\M$ is convex along smooth geodesics, but, again, in applications there is a need to consider also generalized geodesics. It was proved by Chen, \cite{3Chen}, that the K-energy can be defined along any generalized geodesic of class $C^{(1,1)}$, which is  crucial since any two points in the space can be connected by such geodesics. (For this, he rewrites the  definition of $\M$  since the original definition involves four derivatives of $\phi$.) Chen also conjectured that $\M$ would be convex along any generalized geodesic of class $C^{(1,1)}$. This was proved in a joint paper with Robert Berman, \cite{Berman-Berndtsson}.
\begin{thm}\label{Jams} The Mabuchi K-energy is convex along any generalized geodesic of class $C^{(1,1)}$.
\end{thm}
(An alternative proof was later given in \cite{Chen-Long-Paun}.)
Using this we proved the uniqueness of metrics of constant scalar curvature up to flows of holomorphic vector fields:
\begin{thm}\label{Jams2} Let $Z$ be a compact manifold and $L\to Z$ a positive line bundle. Let $\omega_0$ and $\omega_1$ be two K\"ahler metrics on of constant scalar curvature in $c_1[L]$. Then there is a holomorphic vector field on $Z$, with time 1 flow $F$, such that $F^*(\omega_1)=\omega_0$.
\end{thm}
Uniqueness was proved earlier, in case $Z$ has discrete automorphism group in \cite{4Donaldson}. Our result is in fact more general; it treats not only metrics of constant scalar curvature, but general 'extremal metrics', and does not assume that the cohomology class of the metrics is integral. For this, and a discussion of previous work, we refer to \cite{Berman-Berndtsson}.

\section{ Variation of complex structure}
In this section we will mainly consider families of canonically polarized manifolds. We assume that $p:X\to U$ is a smooth proper fibration, that $U$ is a domain in $\C$, and that the fibers $X_t$ have positive canonical bundle. In this setting, $X$ is automatically K\"ahler, since we may find a smooth metric $e^{-\psi}$ on $K_{X/U}$ which is positively curved on fibers, and take a K\"ahler form as $i\ddbar\psi +p^*(\omega)$ where $\omega$ is sufficiently positive on the base. Hence the results from Section \ref{main} apply.

By the Aubin-Yau theorem, \cite{Yau}, \cite{Aubin}, each fiber has a unique K\"ahler-Einstein metric (now with  Ricci curvature equal to -1). This metric can be written
$$
\omega=i\ddbar\phi,
$$
where $e^{-\phi}$ is a metric  on the canonical bundle. Then $e^\phi$ is a metric on the anticanonical bundle, which as in the previous section can be identified with a volume form, and the metric is unique if we normalize so that
$$
e^\phi=(\omega_\phi)^n/n!.
$$
Applying this to each fiber  $X_t$ we get a metric on the relative canonical bundle $K_{X/U}$ that we also denote by $e^{-\phi}$. An important theorem of Schumacher, \cite{Schumacher} implies that $i\ddbar\phi$ is semipositively curved not only along the fibers, but on the total space:
\begin{thm}\label{KESch}   If $e^{-\phi}$ is the normalized K\"ahler-Einstein potential described above, and $c(\phi)$ is defined as in section \ref{main}, then
  $$
  \Box c(\phi)+c(\phi)=|\kappa_\phi|^2
  $$
  on each fiber. As a consequence, $c(\phi)$ is semipositive and strictly positive on each fiber where $\kappa_\phi$ does not vanish identically.
\end{thm}
Here $\Box=\dbar^*\dbar$ is the Laplace operator on functions on a fiber, for the K\"ahler metric $\omega_\phi|_{X_t}$, and $\kappa_\phi$, as defined as in section \ref{main} turns out to be the harmonic representative of the Kodaira-Spencer class. The non negativity part of the statement  follows immediately from the differential equation, since $\Box c(\phi)\leq 0$ at a minimum point, and the strict positivity part is also a well known property of elliptic equations. Since $i\ddbar\phi>0$ along the fibers, the positivity of $c(\phi)$ implies that $i\ddbar\phi$ is positive on the total space.

This result can be seen as an analog of Theorem \ref{Duke}, when $L$ is trivial or a power of the relative canonical bundle. Indeed, that theorem says the Bergman kernel defines a semipositive metric on the relative canonical bundle, at least if it is not identically zero. Tsuji, \cite{Tsuji}, independently proved the semipositivity part of Theorem \ref{KESch}, using Theorem \ref{Duke} and iteration:
\begin{thm}\label{Tsuji} Let $X$ be a family of canonically polarized manifolds. For   any sufficiently positive line bundle $(L,e^{-\psi)}$ on $X$, let
  $e^{-s(\psi)}$ be the Bergman kernel metric on $K_{X/U}+L$.  Define iteratively a sequence of metrics $h_m= e^{-s^m(\psi)}$ on $mK_{X/U}+L$ in this way. Then an appropriate
  renormalization of $h_m^{1/m}$ tends to the metric on $K_{X/U}$ defined by the K\"ahler-Einstein potentials.
  \end{thm} 
Since all the metrics in the iteration have semipositive curvature by Theorem \ref{Duke}, this gives a different proof of the semipositivity part of Theorem \ref{KESch}. Tsuji also applied these arguments to situations when we assume much less positivity along the fibers. 

Both Schumacher's theorem and Theorem \ref{Duke} give some positivity of $K_{X/Y}$ under assumptions of fiberwise positivity. Schumacher's theorem was generalized by Mihai Paun, \cite{Mihai-Alb}, to twisted relative canonical bundles, and even to general adjoint classes, not necessarily integral. He also applied this result to solve a long standing conjecture about the surjectivity of the Albanese map for compact K\"ahler manifolds with nef anticanonical bundle. 
\begin{thm}\label{MihaiAlb} Let $p:X\to Y$ be a holomorphic surjective map between compact K\"ahler manifolds. Let $\beta$ be a semipositive closed $(1,1)$-form on on $X$ and assume that the cohomology class $c_1[K_{X/Y}] +[\beta]$ contains a $(1,1)$-form $\Omega$ which is strictly positive on all fibers $X_y$, for $y$ outside a proper analytic subset, $Y_0$,  of $Y$. Then $c_1[K_{X/Y}]+\beta$ contains a closed semipositive current, which is smooth outside $p^{-1}(Y_0)$.
  \end{thm}

Let us now consider the particular case of a fibration $p:X\to U$ where the fibers are Riemann surfaces of a certain genus $g\geq 2$.
If we choose $L=K_{X/U}$, our vector bundle $E$ with fibers
$$
E_t=H^{1,0}(X_t, K_{X/Y})=H^0(X_t, 2K_{X/U})
$$
is the dual of the bundle with fibers
$$
E^*_t=H^{0,1}(X_t, T^{1,0}(X_t)),
$$
which is the space of infinitesimal deformations of the complex structure on $X_t$. Any positively curved metric on $E$ therefore induces a negatively curved metric on the space of deformations (along $U$). In particular, taking the $L^2$-metric on $E$, defined by the metric on $K_{X/U}$ given by the K\"ahler-Einstein potentials, $e^{-\phi}$, we get the Weil-Peterson metric on $U$. By Schumacher's theorem, $e^{-\phi}$ is (semi)positively curved on $X$, so we get the classical fact that the Weil-Peterson metric has seminegative curvature, and negative curvature where the Kodaira-Spencer class does not vanish, \cite{Ahlfors}, \cite{Royden}, \cite{Wolpert}. Moreover, from Theorem \ref{5pos}, we get an explicit formula for the curvature. Using Schumacher's theorem again (now the differential equation for $c(\phi)$), we can rewrite the formula so that it coincides with the formula found by Wolpert, \cite{Wolpert}, see \cite{3Berndtsson}. The same argument works when $U$ has higher dimension, and then Theorem \ref{Ann2} shows that  the Weil-Peterson metric has dual Nakano-negative curvature, cf \cite{Yau et al}. This is stronger than negative bisectional curvature, which corresponds to Griffiths negativity.

The case of families of higher dimensional canonically polarized manifolds is significantly more complicated. We are then primarily interested in vector bundles $\H^{1,n-1}$ with fibers
$$H^{1,n-1}(X_t, K_{X/U})
$$
and their positivity properties. The reason for this is that the dual of   $H^{1,n-1}(X_t, K_{X/U})= H^{n, n-1}(X_t, (T^*)^{1,0}(X_t))$ is  $H^{0,1}(X_t, T^{1,0})$, the space of infinitesimal deformations of the complex structure on $X_t$. The Kodaira-Spencer map sends the tangent bundle of $U$  into this space, and therefore any metric on $\H^{1,n-1}$ induces a metric on the base. The   Weil-Peterson metric on the base (which can be  obtained this way) was first studied by Siu, \cite{Siu}, who found an explicit formula for its curvature. This was generalized by Schumacher, \cite{Schumacher}, who found a curvature formula 
for the bundle $\H^{1,n-1}$ and also for the other higher direct image bundles $\H^{p,q}$, with $p+q=n$. Both Siu's and Schumacher's formulas contain terms that give an apparent positive contribution to the curvature of the Weil-Peterson metric, but as shown in \cite{Schumacher} and \cite{To-Yeung}, the metrics on all the higher direct images bundles can be combined to give a Finsler metric  of strictly negative sectional curvature that can partly substitute for the Weil-Peterson metric.

All these works focus on the relative canonical bundle and the metric on it given by the K\"ahler-Einstein potential, and in the  case of one dimensional fibers they are special cases   Theorem \ref{5pos}. The generalization to general line bundles $L$ with metrics of fiberwise positive curvature was obtained in \cite{Naumann} and \cite{BMW}. In \cite{BMW}, the formalism was also extended to families of Calabi-Yau manifolds (where related results have also been announced by To and Yeung). Finally, we also mention that, building on work by Lu, \cite{Lu}, Wang,\cite{Wang}, has found a different approach to these problems, which so far works for Calabi-Yau families and then produces a Hermitian version of the Weil-Peterson metric, with negative bisectional curvature. The main idea is to embed the tangent bundle of $U$, not in the dual of $\H^{1,n-1}$ as above, but in the endomorphism bundle of the sum of all the $\H^{p,q}$. 

\section{Positivity of direct image sheaves}

We shall finally discuss extensions of Theorem \ref{Ann2} under less restrictive assumptions. Let $p:X\to Y$ be a surjective holomorphic map between two projective varieties, and let $L\to X$ be a holomorphic line bundle equipped with a singular metric $e^{-\phi}$ with semipositive curvature current $i\ddbar\phi\geq 0$. The 'multiplier ideal sheaf' $\I(\phi)$ is the sheaf of holomorphic functions on $X$ that  are square integrable against $e^{-\phi}$. 
The first problem is that under these general circumstances we no longer get a vector bundle on the base. The role of the vector bundle $E$ is instead played by the (zeroth) direct image sheaf,
$$
\E:=p_*((K_{X/Y}+L)\otimes \I(\phi)).
$$
 
$\E$ is a sheaf over $Y$ whose sections over an  open set $U$ in the base are the sections of $(K_{X/Y}+L)\otimes\I(\phi)$  over $p^{-1}(U)$. 
By a classical theorem of Grauert, $\E$ is  coherent and  it is also torsion free. The consequence of this that we will use is that it is  locally free outside of a subvariety, $Y_0$,  of $Y$ of codimension at least two. This means that outside of $Y_0$, $\E$ is the sheaf of sections of a vector bundle $E$. In the setting of Theorem  \ref{Ann2}, $Y_0$ is empty and $\E$ coincides with the sheaf of sections of $E$ as defined there. In general, let $Y_1$ be a proper subvariety of $Y$ such that $p$ is a submersion outside $p^{-1}(Y_1)$. 
We can now  define a $L^2$-metric  on $E$ over $p^{-1}(Y\setminus (Y_0\cup Y_1))$ as we did when $e^{-\phi}$ was assumed to be smooth. The new feature that appears is that this is now a singular metric. For singular metrics on a vector bundle it seems impossible to define positivity and negativity in terms of a curvature current (cf \cite{Raufi}), but we can circumvent this problem in a way similar to what we did in Theorem \ref{Gent}: We say that a singular metric is negatively curved if the logarithm of the norm of any holomorphic section is plurisubharmonic, and it is positively curved if the dual is negative (cf \cite{Berndtsson-Paun}).

Paun and Takayama extended the notion of singular metrics on vector bundles to the setting of coherent, torsion free sheaves. The definition is simply that a singular metric on such a sheaf is a singular metric on the vector bundle defined by the sheaf outside of $Y_0$. Since negativity (and positivity) of the curvature is defined in terms of plurisubharmonic functions, and since plurisubharmonic functions extend over varieties of codimension at least 2, it turns out that this is a useful definition. Given all this, we have the following theorem of Paun-Takayama, \cite{Paun-Takayama}, and Hacon-Popa-Schell, \cite{Hacon-Popa-Schnell}, which seems to be the most general theorem on (metric) positivity of direct images. 
\begin{thm}\label{HPS} The $L^2$-metric on $E=p_*((K_{X/Y}+L)\otimes \I(\phi))$ over the complement of $Y_0\cup Y_1$ extends to a (singular) metric of nonnegative curvature on the coherent and torsion free sheaf $\E$ .
\end{thm}
(Paun and Takayama proved this when $\I(\phi)$ is trivial on a generic fiber.)
The next theorem is also due to Paun and Takayama:
\begin{thm}\label{PaunTaka}  Any coherent torsion free sheaf which has a (singular) metric of nonnegative curvature is weakly positive in the sense of Viehweg, \cite{Viehweg}. Hence, by Theorem \ref{HPS}, 
  $$
  \E:=p_*((K_{X/Y}+L)\otimes \I(\phi))
  $$
  is weakly positive. 
\end{thm} 
In very general terms, the notion of 'positivity' in algebraic geometry, is often given in terms of the existence of certain holomorphic or algebraic objects. Thus,  e. g. the positivity of a line bundle means that a high power of it has enough sections to give an embedding into projective space. Similarily, Viehweg's weak positivity of a sheaf $\F$ means that certain associated sheaves are generated by global sections. From the analytic point of view, the existence of such sections  should be  a consequence of the metric positivity of curvature, and the previous theorem gives a very general version of this.
We will not go into more details on these matters or give the exact definition of positivity in the sense of Viehweg -- it would require another article (and another author.)

Instead we end by a major application of these results, due to Cao and Paun, \cite{Cao-Paun}.
\begin{thm}\label{CaoPaun} Let $p:X\to A$ be a surjective holomorphic map, where $X$ is smooth projective and $A$ is an Abelian variety. Then we have the inequality for the Kodaira dimensions of $X$ and a generic fiber, $F$,
  $$
  \kappa(X)\geq\kappa(F).
  $$
\end{thm}
This proves the case of the famous 'Itaka conjecture', $\kappa(X)\geq \kappa(F)+\kappa(Y)$, when  $Y$ is an Abelian variety (and hence has Kodaira dimension zero). 
  For a simplification of the proof and a generalization of the result we refer to \cite{Hacon-Popa-Schnell}, which is also a beautiful survey of the field.

 \def\listing#1#2#3{{\sc #1}:\ {\it #2}, \ #3.}


\begin{thebibliography}{9999}

\bibitem{Ahlfors}\listing{Ahlfors, L.}{Curvature properties of Teichm\"uller space}{ J. Analyse Math. 9 (1961), pp. 161-176}   

\bibitem{Aubin}\listing{ Aubin, T,}{ Nonlinear equations on manifolds. Monge Amp\`ere equations}{ Grundlehren der Mathematischen Wissenschaften 252, Springer Verlag, (1982)}

  
\bibitem{Bando-Mabuchi}\listing{Bando, S. and  Mabuchi, T.}{ Uniqueness of Einstein K\"ahler
metrics modulo connected group actions.}{ Algebraic geometry, Sendai,
1985, 11\textendash{}40, Adv. Stud. Pure M ath., 10, North-Holland,
Amsterdam, 1987.}
\bibitem{Bedford-Taylor}\listing{Bedford, E. and Taylor, B.A.}{  The Dirichlet problem for the complex Monge-Amp\`ere equation}{Invent. Math. 103 (1976), pp. 1-44} 

\bibitem{Berman et al}\listing{Berman R.J: Eyssidieux, P: Boucksom, S; Guedj,
V; Zeriahi, A}{ K\"ahler-Einstein metrics and the K\"ahler-Ricci flow on
log Fano varieties}{ J. Reine Angew. Math., to appear}

\bibitem{Berman}\listing{Berman, R.J.}{Sharp inequalities for determinants of Toeplitz operators and dbar-Laplacians on line bundles}{arXiv:0905.4263}

\bibitem{2Berman}\listing{Berman, R. J.}{K-polystability of Q-Fano varieties admitting K\"ahler-Einstein metrics.}{ Invent. Math. 203 (2016), no. 3, 973-1025.}

  \bibitem{Berman-Berndtsson}\listing{Berman, R. J. and Berndtsson B.}{ Convexity of the K-energy on the space of K\"ahler metrics and uniqueness of extremal metrics}{ J. Amer. Math. Soc. 30 (2017), pp. 1165-1196.}
 
\bibitem{1Berndtsson}\listing{Berndtsson, B}{Curvature of vector bundles
    associated to holomorphic fibrations  }{Ann. of Math.  169 2009, pp
    531-560 }
\bibitem{2Berndtsson}\listing{Berndtsson B}{Subharmonicity of the Bergman kernel and some other functions associated to pseudoconvex domains}{Ann Inst Fourier, 56 (2006) pp 1633-1662}
\bibitem{3Berndtsson}\listing{Berndtsson, B}{ Strict and nonstrict positivity of direct image bundles. }{ Math. Z. 269 (2011), no. 3-4, 1201-1218.}
\bibitem{4Berndtsson}\listing{Berndtsson, B.}{A Brunn-Minkowski inequality for Fano manifolds and some uniqueness theorems in K\"ahler geometry}{ Invent. Math. 200 (2015), pp. 149-200}


  
\bibitem{Berndtsson-Lempert}\listing{Berndtsson, B and Lempert, L}{A proof of the Ohsawa-Takegoshi theorem with sharp estimates }{ J. Math. Soc. Japan 68 (2016), pp. 1461-1472}

\bibitem{Berndtsson-Paun}\listing{Berndtsson, B and Paun, M}{Bergman kernels and the pseudoeffectivity of relative canonical bundles.}{ Duke  Math. J. 145 ,2008 no. 2, 341-378.}

\bibitem{BMW}\listing{Berndtsson, B., Paun, M. and Wang, X.}{Algebraic fiber spaces and curvature of higher direct images}{ Arxiv 1704.02279}



\bibitem{Blocki}\listing{B{\l}ocki, Z.}{Suita's conjecture and the Ohsawa-Takegoshi extension theorem}{ Invent. Math. 193 (2103), pp. 149-158}  


\bibitem{2Blocki}\listing{B{\l}ocki, Z.}{Bergman kernel and pluripotential theory}{Analysis, complex geometry, and mathematical physics: in honor of Duong H. Phong,pp 1-10,
Contemp. Math., 644, Amer. Math. Soc., Providence, RI, 2015. }Algebraic fiber spaces and curvature of higher direct images


\bibitem{3Blocki}\listing{B{\l}ocki, Z.}{ On geodesics in the space of K\"ahler metrics,}{
Proceedings of the \textquotedbl{}Conference in Geometry\textquotedbl{}
dedicated to Shing-Tung Yau (Warsaw, April 2009), in \textquotedbl{}Advances
in Geometric Analysis\textquotedbl{}, ed. S. Janeczko, J. Li, D. Phong,
Advanced Lectures in Mathematics 21, pp. 3-20, International Press,
2012 }


  
 \bibitem{Brascamp-Lieb}\listing  {H J Brascamp and E H  Lieb} {On
  extensions of the Brunn-Minkowski and Pr\'ekopa-Leindler theorems,
  including inequalities for log concave functions, and with an
  application to the diffusion equation.}  {J. Functional Analysis  22
   (1976), no. 4, 366--389}
 \bibitem{Cao-Paun}\listing{ Cao, J. and Paun, M.}{ Kodaira dimension of algebraic fiber spaces over abelian varieties.}{Invent. Math. 207 (2017), pp. 345-387}   
\bibitem{Chen}\listing{Chen, X.X}{ The space of K\"ahler metrics }{ J. Diff.
Geom. 56 (2000), 189-234}
\bibitem{3Chen}\listing{Chen, X.X}{ On the lower bound of the Mabuchi energy and
  its application}{ Int. Math. Res. Not. 2000, no. 12, 607-623}
  \bibitem{CDS}\listing{ Chen, X.X., Donaldson, S.K. and Song, S.}{ K\"ahler-Einstein metrics and stability}{ Int. Math. Res. Not. 2014, pp. 2119-2125}

 \bibitem{Chen-Long-Paun} \listing{Chen, XX, Li, L and Paun M.}{ Approximation of weak geodesics and subharmonicity of Mabuchi energy}{ Ann. Fac. Toulouse Math. (6) 25 (2016), pp. 935-957}
  
 \bibitem{Dario}\listing{ Cordero-Erausquin, D.} {On Berndtsson's generalization of Pr\'ekopa's theorem.}{ Math. Z. 249 (2005), no. 2, 401-410.}
 \bibitem{Dario-Boaz}\listing{Cordero-Erausquin, D. and Klartag, B.}{Interpolations, convexity and geometric inequalities}{Geometric Aspects of Functional Analysis, SLN 2050, Springer, Heidelberg, 2012}  

\bibitem{Darvas-Lempert}\listing{Darvas, T; Lempert, L.}{ Weak geodesics in the space
of K\"ahler metrics}{ Math. Res. Lett. 19 (2012), no. 5, 1127\textendash{}1135.}       
\bibitem{Demaillynotes}\listing{Demailly, J-P}{Complex analytic and differential geometry.}{Notes available on Demailly's webpage}

\bibitem{Ding}\listing{Ding, W-Y.}{Remarks on the existence problem for positive K\"ahler-Einstein metrics}{ Math. Ann. 282 (1988), pp. 463-471}
  
\bibitem{Donaldson}\listing{Donaldson, S. K}{ Symmetric spaces, K\"ahler geometry
and Hamiltonian dynamics}{ In Northern California Symplectic Geometry
Seminar , volume 196 of Amer. Math. Soc. Transl. Ser. 2 , pages 13\textendash{}33.}
Amer. Math. Soc., Providence, RI, 1999
\bibitem{2Donaldson}\listing{Donaldson, S. K.}{ Scalar curvature and projective embeddings II.}{
  Q. J. Math. 56 (2005), no. 3, 345--356.}
  \bibitem{3Donaldson}\listing{Donaldson, S. K.}{ K\"ahler metrics with cone singularities along a divisor}{ Essys in mathematics and its applications, pp. 49-7949-79, Springer, Heidelberg 2012}

  
\bibitem{4Donaldson}\listing{Donaldson, S.K}{ Scalar curvature and projective embeddings,
I}{ J. Diff. Geom. 59 (2001), 479\textendash{}522}
\bibitem{Dubuc}\listing{ Dubuc, S.}{ Criteres de convexite et inegalites integrales.}{ Ann. Inst. Fourier Grenoble 27 (1) (1977), 135-165.}







  
\bibitem{Fujita}\listing{ Fujita, T.}{On K\"ahler fiber spaces over
  curves. }{ J. Math. Soc. Japan  30  (1978), no. 4, 779--794}
\bibitem{Gardner}\listing{Gardner, R. J.}{ The Brunn-Minkowski
    Inequality}{BAMS, 39 (2002), pp 355-405}




  
\bibitem{Griffiths}\listing
{Griffiths, P A}{Curvature properties of the Hodge
    bundles (Notes written by Loring Tu) }{ Topics in Transcendental Algebraic Geometry, Annals of
    Mathematics Studies, Princeton University Press 1984}
  

      \bibitem{Guan-Zhou}\listing{Guan, Q. and Zhou, X.}{A solution of an $L^2$-extension problem with an optimal estimate and applications }{Ann. of Math. 181 (2015), pp. 1139-1208}

      \bibitem{Hacon-Popa-Schnell}\listing{Hacon, C., Popa, M. and Schnell, C.}{Algebraic fiber spaces over analytic varieties: around a recent theorem by Cao and Paun}{Arxiv 1611.08768}
\bibitem{Hormander}\listing{H\"ormander, L.}{ $L^2$-estimates and existence theorems for  the $\dbar$-operator}{Acta Math. 113 (1965), pp. 89-152}
        
\bibitem{Kawamata}\listing{ Kawamata, Y.}{ Subadjunction of log canonical divisors}{Amer. J. Math. 120 (1998), pp. 893-899.}        
\bibitem{Lempert-Vivas}\listing{Lempert, L; Vivas, L}{ Geodesics in the space of K\"ahler
metrics.}{ Duke Math. J. Volume 162 , Number 7 (2013), 1369-138}


\bibitem{Kiselman}\listing{Kiselman, C. }{ The partial Legendre transform for plurisubharmonic functions}{ Invent. Math. 49 (1978), pp. 137-148}

\bibitem{Yau et al}\listing{ Liu,K;  Sun,X and   Yau S-T}{Recent Development
    on the Geometry of the Teichmuller and Moduli Spaces of Riemann
    Surfaces and Polarized Calabi-Yau Manifolds  }{Surveys in differential geometry. Vol. XIV. Geometry of Riemann surfaces and their moduli spaces, 221-259, Surv. Differ. Geom., 14, Int. Press, Somerville, MA, 2009.} 


 \bibitem{Lu}\listing{ Lu, Z.}{On the Hodge metric of the universal deformation space of Calabi-Yau three-folds}{ Jour. Geom. Anal. 11 (2001), pp. 103-118} 



  
\bibitem{Mabuchi}\listing{Mabuchi, T}{ Some symplectic geometry on compact K\"ahler
  manifolds}{ I. Osaka J. Math. , 24(2):227\textendash{}252, 1987.}

 

\bibitem{2Mabuchi}\listing{Mabuchi, T.}{ K -energy maps integrating Futaki invariants.}{
Tohoku Math. J. (2) 38 (1986), no. 4, 575-593}





\bibitem{Maitani-Yamaguchi}\listing{ Maitani, F and Yamaguchi, H}{ Variation of Bergman metrics on Riemann surfaces}{ Math Ann 330 (2004) pp 477-489}

\bibitem{Mihai-Alb}\listing{Paun, M.}{Relative adjoint transcendental classes and Albanese map of compact K\"ahler manifolds with nef Ricci curvature}{Advanced Studies in Pure Mathematics, Math. Soc of Japan, to appear}

\bibitem{Moser}\listing{Moser, J.}{ A sharp form of an inequality by N. Trudinger}{ Indiana Univ. Math. J. 20 (1970/71)}

\bibitem{Naumann}\listing{Naumann, P.}{Curvature of higher direct images}{ Arxiv 1611.09117}

\bibitem{Ohsawa-Takegoshi}\listing{ Ohsawa, T. and Takegoshi, K.}{On the extension of $L^2$-holomorphic functions}{ Math Z. 195 (1987), pp. 197-204}  

  \bibitem{Onofri}\listing{Onofri, E.}{ On the positivity of the action in a theory of random surfaces}{ Comm. Math. Phys. 86 (1982), pp. 321-326}
   
  \bibitem{Paun-Takayama}\listing{Paun, M and Takayama, S}{ Positivity of relative pluricanonical bundles and their images}{arXiv:1409.5504}

 \bibitem{Phong-Sturm}\listing{Phong, D. H. and Sturm, J.}{Lectures on stability and constant scalar curvature}{ Current developments in mathematics, 2007, pp. 101-176 (Int. Press, Somerville, MA 2009)}   

\bibitem{Prekopa}\listing{Pr\'ekopa, A.}{ On logarithmic concave measures and functions.}{ Acta Sci. Math. (Szeged) 34 (1973), pp. 335-343}       
        
\bibitem{Raufi}\listing{Raufi, H.}{Singular metrics on holomorphic vetor bundles.}{Ark. Mat. 53 (2015), pp. 359-382}

  \bibitem{2Raufi}\listing{Raufi, H.}{ Log concavity for matrix-valued functions and a matrix-valued Pr\'ekopa theorem}{Arxiv 1311.7343}

      \bibitem{David-Ross}\listing{ Ross, J. and Witt-Nystr\"om, D.}{ Harmonic Discs of Solutions to the Complex Homogeneous Monge-Amp\`ere Equation}{ Publ. Math. Inst. Hautes \'Etudes Sci. 122 (2015), pp. 315-335. }

 \bibitem{2David-Ross}\listing{    Ross, J. and Witt-Nystr\"om, D.}{ On the maximal rank problem for the complex homogeneous Monge-Amp\`ere equation}{ Axiv, 1610.02280 }     

\bibitem{Royden}\listing{Royden, H.}{ Intrinsic metrics on Teichm\"uller space}{ Proc. Int. Cong. Math. 2 (1974), pp. 217-221}        

        

\bibitem{Schumacher}\listing{ Schumacher, G}{Positivity of relative
    canonical bundles and applications}{Invent. Math. 190 (2012), pp 1-56}
  


        

\bibitem{Semmes}\listing{ Semmes, S}{ Interpolation of Banach spaces, differential geometry and differential equations.
}{Rev. Mat. Iberoamericana 4 (1988), no. 1, 155--176}
        

\bibitem{Siu}\listing{Siu, Y-T}{Curvature of the Weil-Petersson metric
  in the moduli space of compact K\"ahler-Einstein manifolds of negative
first Chern class}{In: Contributions to several complex variables, ed
  A Howard and P-M Wong, Vieweg 1986}

\bibitem{To-Yeung}\listing{To, W-K. and Yeung, S-K.}{Finsler metrics and Kobayashi hyperbolicity of the moduli space of canonically polarized manifolds}{ Ann. of Math. 181 (2015), pp. 547-586}


  
\bibitem{Tsuji}\listing{Tsuji, H.}{   Dynamical construction of K\"ahler-Einstein metrics.}{ Nagoya Math. J. 199 (2010), 107-122.} 


\bibitem{Viehweg}\listing{Viehweg, E.}{ Quasi-projective moduli for polarized manifolds}{ Ergebnisse der Mathematik und ihre grenzgebiete (3) 30, Springer-Verlag Berlin 1995}

  
\bibitem{1Wang}\listing{Wang, X.}{A curvature formula associated to a family of pseudoconvex domains}{ Ann. Inst. Fourier (Grenoble) 67 (2017), pp 269-313}


\bibitem{Wang}\listing{Wang, X.}{ Curvature restrictions on a manifold with a flat Higgs bundle}{ Arxiv 1608.00777}







  
\bibitem{Wolpert}\listing{Wolpert, S}{Chern forms and the Riemann tensor for the moduli space of curves}{Invent. Math 85 1986, pp. 119-145 }
\bibitem{Yau}\listing{Yau, S.T.}{ On the Ricci curvature of a compact K\"ahler manifold and the complex Monge-Amp\`ere equation.}{ Comm. Pure Applied Math. 31 (1978), pp. 339-411}

\end{thebibliography}
\end{document}